\documentclass[a4paper,12pt,reqno]{amsart}
\usepackage{amsmath,amssymb,latexsym,amsfonts,amsthm,comment, xspace, enumerate} 
\usepackage[english]{babel}
\DeclareMathAlphabet\mathbfcal{OMS}{cmsy}{b}{n}
\input amssym.def
\input amssym
\usepackage{latexsym}
\usepackage{amsbsy}
\usepackage{amssymb}
\usepackage{amsfonts}
\usepackage{amsmath}
\usepackage[latin1]{inputenc}
\usepackage{graphicx}
\usepackage{a4wide}
\usepackage{mathtools}
\usepackage{tikz}
\usepackage{amsthm}
\usepackage[english]{babel}

 \newtheorem{thm}{Theorem}[section]
 \newtheorem{defn}[thm]{Definition}
 \newtheorem{lm}[thm]{Lemma}
 
 \newtheorem{cor}[thm]{Corollary}

 \newcommand{\beq}{\begin{equation}}
 \newcommand{\eeq}{\end{equation} }

 \newcommand{\bea}{\begin{eqnarray}}
 \newcommand{\eea}{\end{eqnarray}}

 \newcommand{\beas}{\begin{eqnarray*}}
 \newcommand{\eeas}{\end{eqnarray*}}

 \newcommand{\beqs}{\begin{equation*}}
 \newcommand{\eeqs}{\end{equation*}}

 \newcommand{\bi}{\begin{itemize}}
 \newcommand{\ei}{\end{itemize}}

 \newcommand{\ben}{\begin{enumerate}}
 \newcommand{\een}{\end{enumerate}}

 \newcommand{\ba}{\begin{array}}
 \newcommand{\ea}{\end{array}}

\newcommand{\NN}{\mathbb N}
\newcommand{\CC}{\mathbb C}
\newcommand{\RR}{\mathbb R}
\newcommand{\ZZ}{\mathbb Z}

\newcommand{\SSS}{\mathcal S}

\begin{document}

\title{On the iterates of the Laguerre operator}

\author{Smiljana Jak\v si\'c}
           \address{Faculty of Forestry, University of Belgrade, Kneza Vi\v seslava 1, Belgrade, Serbia}
 \email{smiljana.jaksic@sfb.bg.ac.rs}

\author{Stevan Pilipovi\'c }
          \address{Department of Mathematics, Faculty of Sciences, University of Novi Sad, Trg Dositeja Obradovi\'ca 4, Novi Sad, Serbia}
          \email{stevan.pilipovic@dmi.uns.ac.rs}

   \author{Nenad Teofanov}
  \address{Department of Mathematics, Faculty of Sciences, University of Novi Sad, Trg Dositeja Obradovi\'ca 4, Novi Sad, Serbia}
  \email{tnenad@dmi.uns.ac.rs@dmi.uns.ac.rs}

\author[\DJ. Vuckovi\' c]{\DJ or\dj e Vu\v{c}kovi\' c}
\address{\DJ or\dj e Vu\v ckovi\'c, Technical Faculty ``Mihajlo Pupin'', \DJ ure \DJ akovi\'{c}a bb, 23000 Zrenjanin, Serbia}
\email{djordjeplusja@gmail.com}

\begin{abstract}
We use the iterates of the Laguerre operator to introduce Pilipovi\'c spaces on positive orthants. It is shown that such spaces coincide with $G-$type spaces $g_\alpha^\alpha(\RR^d_+)$ and $G_\alpha^\alpha(\RR^d_+)$, when $\alpha > 1$, and  $\alpha \geq  1$, respectively. However, in contrast to  $G$-type spaces, Pilipovi\'c spaces are nontrivial below the critical index $\alpha = 1$. We also remark that there is a natural isomorphism between subspaces of Pilipovi\'c spaces on $\RR ^d$ consisting of even functions, and  Pilipovi\'c spaces on positive orthants.
\end{abstract}



\maketitle

\section{Introduction} \label{sec:1}

Counterparts of Schwartz functions and of Gelfand-Shilov type functions on $(0,\infty)$ were studied in \cite{D0,D1}. These results were thereafter generalized in \cite{SSB}, were the  test spaces of ultradifferentiable functions of Roumieu and Beurling  type, denoted by $G_\alpha^\alpha(\RR^d_+)$, and $g_\alpha^\alpha(\RR^d_+)$, $\alpha\geq 1$, respectively, are introduced and characterized via decay properties of their coefficients in Laguerre expansions.

Toft in \cite{Toft1} studied Pilipovi\'{c} spaces $\mathbfcal S_{\alpha}^{\alpha}(\mathbb R^d)$, $ \alpha >0$, by considering the iterates of the linear Harmonic oscillator. These spaces are non trivial even if $0<\alpha<\frac{1}{2}$, whereas when $\alpha \geq \frac{1}{2}$ they coincide with the Gelfand-Shilov spaces usually denoted by $\mathcal S_{\alpha}^{\alpha}(\mathbb R^d)$, see e.g. \cite{GS2}.
Recently, Toft, Bhimani and Manna demonstrated that Pilipovi\'c spaces on $\RR^d$, and their distribution spaces naturally appear in the study of well-posednedness of certain Schr\"odinger equations, which are ill-posed in the class of Gelfand-Shilov type distributions, see \cite{TBM}.

In this report we use the iterates of the Laguerre operator to introduce Pilipovi\'c spaces on positive orthants, and relate them with with
$G-$type spaces $g_\alpha^\alpha(\RR^d_+)$ and $G_\alpha^\alpha(\RR^d_+)$ of ultradifferentiable functions when $\alpha > 1$, and  $\alpha \geq  1$, respectively.

The proofs will be given in our forthcoming contribution \cite{JPTV}, where we also study natural isomorphisms between Pilipovi\'c spaces on $\RR^d$ consisting of even functions and Pilipovi\'c spaces on $\RR^d _+$. Related questions in the context of $G$-type spaces were considered in \cite{SSSB}.

\subsection{Notation}

By $\NN$, $\ZZ$, $\RR$ and $\CC$ the sets of positive
integers, integers, real and complex numbers, respectively;
$\NN_0=\NN \cup\{0\}$,  $\NN^d_0= (\NN \cup\{0\})^d$,
$\RR_+=(0,\infty)$, $\RR^d_+=(0,\infty)^d$ and $\overline{\RR^d_+}=[0,\infty)^d$.
The standard multi-index
notation will be used: for $\alpha\in\NN_0^d$ and $x\in\RR^d$ (or $x\in\overline{\RR^d_+}$), $|\alpha|:=\sum_{i=1}^d \alpha_j$, $x^\alpha=x_1^{\alpha_1}...x_d^{\alpha_d}$, and $\partial^\alpha=\partial^{\alpha_1}/{\partial x_1^{\alpha_1}}...\partial^{\alpha_d}/\partial x_d^{\alpha_d}$. If  $x,\gamma\in\overline{\RR^d_+}$, then $x^{\gamma}= x_1^{\gamma_1}...x_d^{\gamma_d}$, and if
$x_j=0$ and $\gamma_j=0$, $j=1,\dots,d$, we use the convention $0^0=1$.

We use the standard notation for common spaces of functions, distributions, and sequences.
For example, $ \ell_p $, $ p\geq 1$, denotes the Banach space of Lebesgue sequences, endowed with the  usual norm:
$$
\{a_n \}_{n \in \NN ^d _0} \in \ell_p \qquad \Longleftrightarrow \qquad \| a_n \|_{\ell^p}  = ( \sum_n |a_n|^p )^{1/p} <\infty,
$$
and with the usual modification when $p=\infty.$

A  measurable function belongs to the Hilbert space  $L^2(\RR^d_+)$ if
$$ \|\cdot\|_{L^2(\RR^d_+)} := (\int_{\RR^d_+} |f(x)|^2 dx)^{1/2}<\infty,
$$
and the scalar product is given in the usual way.
By $\SSS(\mathbb{R}_+^d)$ we denote the space of all smooth functions on positive orthant, $f\in \mathcal{C}^{\infty}(\RR^d_+)$ such that all of its derivatives $D^pf$, $p\in\NN^d_0$, extend to continuous functions on
$\overline{\RR^d_+}$, and
$$\sup_{x\in\mathbb{R}^d_+}x^k|D^pf(x)|<\infty\;,\qquad  k,p\in\mathbb{N}_0^d.$$
With this system of seminorms $\mathcal S(\mathbb{R}_+^d)$ becomes an $(F)$-space, see e.g. \cite{Sm}.


\section{Preliminaries} \label{sec:2}

\subsection{Sequence spaces}

We use the notation from \cite{Toft1}.

Let $\alpha>0 $,  $h>1$, and put $\vartheta_{h,\alpha}(n)= e^{h|n|^{1/(2\alpha)}}, $  $ n \in \NN^d_0$.  By $\ell^{p}_{[\vartheta_{h,\alpha}]}(\NN_0^d)$, $p \geq 1$,
we denote the space of all complex sequences $\{a_n\}_{n\in\NN^d_0}$
for which
\begin{equation} \label{eq:norm-sequence}
\|\{a_n\}_{n\in\mathbb \NN^d_0 }\|_{\ell^{p}_{[\vartheta_{h,\alpha}]}(\NN_0^d)}:=\|\{|a_n|\cdot \vartheta_{h,\alpha}(n)\}_{n\in\NN^d_0}\|_{l^{p}}<\infty.
\end{equation}
Equipped with this norm $\ell^{p}_{[ \vartheta_{h,\alpha}(n)}(\NN_0^d)$ becomes a Banach space, with the usual modification when $p = \infty$.

We also introduce
\begin{equation*}
\ell_{\alpha}(\NN^d_0)=\bigcup_{h>0} \ell^{\infty}_{[\vartheta_{h,\alpha}]}(\NN_0^d), \qquad
  \ell_{0,\alpha}(\NN_0^d) =\bigcap_{h>0} \ell^{\infty}_{[\vartheta_{h,\alpha}]}(\NN_0^d),
\end{equation*}
and endow these locally convex spaces with inductive and projective limit topology, respectively.

For consistency, if $\alpha=0 $, then $\ell_0(\NN^d_0)$ denotes the space of all sequences
$\{a_n\}_{n\in\mathbb \NN^d_0 }$ such that $a_n\not=0$ for only finitely many $n \in \NN_0^d$, and $\ell_{0,0}(\NN^d_0)=\{0\}.$

\par

Observe that different choices of $p \in [1,\infty]$ give rise to equivalent norms in \eqref{eq:norm-sequence}, cf
\cite{Toft1}. In particular, if
$\{a_n\}_{n\in\mathbb \NN^d_0} \in \ell^{\infty}_{[\vartheta_{h_1,\alpha}]}(\NN_0^d) $ then there exists $h_1>0$ and $C_1>0$ such that
$$\|\{a_n\}_{n\in\mathbb \NN^d_0\}}\|_{\ell^{2}_{[\vartheta_{h_1,\alpha}]}(\NN_0^d)}=\Big(\sum_{n\in\NN^d_0}|a_n|^2 \vartheta_{h_1,\alpha}(n)\Big)^{1/2}\leq  C_1 \|\{a_n\cdot \vartheta_{h,\alpha}(n)\}_{n\in\NN^d_0}\|_{l^{\infty}},$$
and conversely,
if $$\{a_n\}_{n\in\mathbb \NN^d_0} \in \ell^{2}_{[\vartheta_{h,\alpha}]}(\NN_0^d)$$ then there exist $C_2,h_2>0$ such that
$$\|\{a_n\}_{n\in\mathbb \NN^d_0\}}\|_{\ell^{\infty}_{[\vartheta_{h_2,\alpha}]}(\NN_0^d)}\leq C_2 \|\{a_n\cdot \vartheta_{h,\alpha}(n)\}_{n\in\NN^d_0}\|_{l^{2}}. $$

\par

\subsection{Laguerre functions and Laguerre operator} \label{lager}

For $j\in\mathbb{N}_0$ and $\gamma\geq 0$, the $j$-th Laguerre
polynomial of order $\gamma$ is defined by
$$L_j^\gamma(x)=\frac{x^{-\gamma}e^x}{j!}\frac{d^j}{dt^j}(e^{-x}x^{\gamma+j}),\qquad x \geq 0.$$
The $j$-th Laguerre function is then given by
$l_j(x)=L_j(x)e^{-x/2}$, $ x \geq 0$.
The $d$-dimensional Laguerre polynomials and Laguerre
functions, are given by tensor products
$L_n(x)=\prod_{l=1}^d L_{n_l}(x_l)$, and
$l_n(x)=\prod_{l=1}^dl_{n_l}(x_l)$, $n\in\NN^d_0$,
respectively.
We are interested in the case $\gamma=0$, and we write $L_n$ and $l_n$ instead of $L_n^0$ and $\l_n^0$, respectively.

\par

Laguerre functions $\{l_n\}_{n\in\NN^d_0}$ form an orthonormal basis in $L^2(\mathbb{R}^d_+)$, cf. e.g. \cite{Thang, Wong}. Thus every $f\in L^2(\mathbb{R}^d_+)$ can be expressed  via its Laguerre series expansion,
and by $a_n (f)$
we denote the corresponding  Laguerre coefficients:
$$a_n(f)=\int_{\mathbb{R}^d_+}f(x)l_{n}(x)dx, \qquad n\in\mathbb{N}^d_0.$$

Similarly as it is done in \cite{Toft1},
we may consider the spaces of formal Laguerre series  expansions
$f=\sum_{n\in\NN_0^d} a_{n} l_n$ that correspond to the sequence $\{a_n \}_{n\in\NN_0^d}$ from $\ell_{\alpha}(\NN_0^d)$ or  $\ell_{0,\alpha}(\NN_0^d)$, $\alpha \geq 0$.
These spaces will be denoted by $\mathcal G_{\alpha}(\RR_{+}^d)$ and $\mathcal G_{0,\alpha}(\RR_{+}^d)$, respectively, and the natural mapping
\begin{multline} \label{T}
T: \ell_{\alpha}(\NN_0^d) \rightarrow\mathcal G_{\alpha}(\RR_{+}^d) \quad
(T: \ell_{0,\alpha}(\NN_0^d) \rightarrow\mathcal G_{0,\alpha}(\RR_{+}^d)),
\\
T(\{a_n \}_{n\in\NN_0^d})=\sum_{n\in\NN_0^d} a_n l_n
\end{multline}
induces  topology on the target space.

Next we introduce the Laguerre operator:
$$
E=-\sum_{j=1}^d \Big(x_j\frac{\partial ^2}{\partial x_j^2}+\frac{\partial}{\partial x_j}-\frac{x}{4}+\frac{1}{2}\Big).
$$
Then the Laguerre functions are the eigenfunctions of $E$:
$$
E ( l_n(x))=\Big(\sum_{i=1}^d  n_j\Big )\cdot l _n(x)=|n| \cdot l _n(x),\qquad x>0,
$$
see \cite[(11), p.188]{Ed}, and for $N\in \mathbb N$, and $n\in \NN_0^d$ we have
\begin{equation}\label{Eigen}
E^{N} ( l_n(x))=|n|^N l_{n}(x).
\end{equation}

\section{Pilipovi\'c spaces on positive orthants}

Before we define spaces of ultradifferentiable functions in terms of  the iterates of the  Laguerre operator
we give a simple but useful lemma.

\begin{lm} \label{lm:S+}
Let $u $ be a smooth and square integrable function on $ \RR_{+}^d$, i.e.
$u\in L^2(\RR_{+}^d) \cap C^{\infty}(\RR^d_{+})$. If $E^Nu\in L^2(\RR_{+}^d)$ for every $N\in\NN$,
then $u\in \mathcal S(\RR_{+}^d).$
\end{lm}

The proof follows from \eqref{Eigen} and the characterization of $\SSS(\mathbb{R}_+^d)$ via coefficients in Laguerre expansions of its elements.

\begin{defn} \label{def:Galphah}
Let $h>0$, and $\alpha>0$. Then a smooth functions $ f $ belongs to the space $ \mathbf G^{\alpha,h}_{\alpha,h}(\RR^d_+)$, if
\begin{equation}\label{eta}
\eta_h^\alpha(f) : =\sup_{N\in\NN}\frac{\Vert E^N f\Vert_{L^2(\RR^d_+)}}{h^{|N|} N!^\alpha}<\infty.
\end{equation}
\end{defn}

From Lemma \ref{lm:S+} it immediately follows that $\mathbf G^{\alpha,h}_{\alpha,h}(\RR^d_+) \subset \mathcal S(\RR^d_{+})$. $\mathbf G^{\alpha,h}_{\alpha,h}(\RR^d_+)$  is a Banach space with the norm $\eta_h^\alpha$, and clearly, if $0<\alpha < \beta$, then $\mathbf G^{\alpha,h}_{\alpha,h}(\RR^d_+)$ is continuously injected into $\mathbf G^{\beta,h}_{\beta,h}(\RR^d_+)$.

The Laguerre functions  $l_p$, $p\in\NN^d_0$, satisfy (\ref{eta}) for every $h>0$ and $\alpha>0$,
since $$\sup_{N\in\NN}\frac{\Vert E^N l_p\Vert_{L^2(\RR^d_+)}}{h^{N} N!^{\alpha}}=\sup_{N\in\NN}\frac{(\frac{|p|}{h})^N}{N!^{\alpha}}<\infty.$$
It follows that  $\mathbf G^{\alpha,h}_{\alpha,h}(\RR^d_+)$ are non-trivial spaces since they contain at least finite linear combinations of the Laguerre functions.

We notice that even if $\alpha=0$, and if $h>1$, then $l_p\in \mathbf G^{0,h}_{0,h}(\RR^d_+)$
for all $p\in\NN^d$ such that $|p|\leq h$.

Now we are ready to introduce Pilipovi\'c spaces on positive orthants.

\begin{defn} \label{def:PilipovicOrthant}
Let $\alpha>0$. Then the Pilipovi\'c spaces on positive orthant $\RR^d_+$, $\mathbf{G}^{\alpha}_{\alpha}(\RR^d_+)$ and  $\mathbf{g}^{\alpha}_{\alpha}(\RR^d_+)$ are defined to be the union, respectively the intersection of $ \mathbf G^{\alpha,h}_{\alpha,h}(\RR^d_+)$, with respect to $h>0$:
\begin{equation*}
 \mathbf G^{\alpha}_{\alpha}(\RR^d_+)=
\bigcup_{h>0}  \mathbf
G^{\alpha,h}_{\alpha,h}(\RR^d_+), \qquad \text{and} \qquad
\mathbf{g}^{\alpha}_{\alpha}(\RR^d_+)=
\bigcap_{h>0}  \mathbf
G^{\alpha,h}_{\alpha,h}(\RR^d_+).
\end{equation*}

\end{defn}

The spaces $\mathbf{G}^{\alpha}_{\alpha}(\RR^d_+)$ and  $\mathbf{g}^{\alpha}_{\alpha}(\RR^d_+)$ are
endowed with  projective and inductive limit topologies, respectively.
Thus, $\mathbf G^{\alpha}_{\alpha}(\RR^d_+)$ and $\mathbf{g}^{\alpha}_{\alpha}(\RR^d_+)$ are $(F)$-spaces continuously injected into $\SSS(\RR^d_+)$, and if $0<\alpha < \beta$ we have
\begin{equation} \label{poredak}
\mathbf{g}^{\alpha}_{\alpha}(\RR^d_+)\hookrightarrow \mathbf{G}^{\alpha}_{\alpha}(\RR^d_+)\hookrightarrow  \mathbf{g}^{\beta}_{\beta} \hookrightarrow \mathbf{G}^{\beta}_{\beta}(\RR^d_+),
\end{equation}
where $\hookrightarrow$ denotes dense and continuous inclusion.

Next we state our main result.

\begin{thm}\label{glavna} Let $\alpha > 0$, $f\in L^2(\mathbb R_+^d) \cap C^\infty(\mathbb R_+^d)$, and let
 $\{a_n(f)\}_{n\in\NN^d_0} $ be the sequence of Laguerre coefficients of $f$. Then the following is true:
\begin{itemize}
\item[(i)] $ f \in \mathbf G^{\alpha}_{\alpha}(\RR^d_+)$  if and only if
 $\{a_n(f)\}_{n\in\NN^d_0}\in \ell_{\alpha/2}(\NN^d_0)$
\item[(ii)] $ f \in \mathbf g^{\alpha}_{\alpha}(\RR^d_+)$  if and only if
 $\{a_n(f)\}_{n\in\NN^d_0}\in \ell_{0,\alpha/2}(\NN^d_0)$.
 \end{itemize}
\end{thm}

The proof will be given in a separate contribution, \cite{JPTV}

By Theorem \ref{glavna} it follows that $\mathbf{G}^{\alpha}_{\alpha}(\RR^d_+)=\mathcal{G}_{\alpha}(\RR^d_+)$ and  $\mathbf{g}^{\alpha}_{\alpha}(\RR^d_+)=\mathcal{G}_{0,\alpha}(\RR^d_+)$ where
$\mathcal{G}_{\alpha}(\RR^d_+)$ and  $\mathcal{G}_{0,\alpha}(\RR^d_+)$ are intriduced in subsection \ref{lager}.
Moreover, the mapping $T:{\ell_{\alpha}(\NN_0^d)}\rightarrow \mathbf G^{\alpha}_{\alpha}(\RR^d_{+})$ (resp. $T:{\ell_{0, \alpha}(\NN_0^d)}\rightarrow \mathbf g^{\alpha}_{\alpha}(\RR^d_{+})$ ) is a well defined topological isomorphism, see \eqref{T} for the definition of $T$.

\subsection{$G$-type spaces as Pilipovi\'c spaces}

Following \cite{SSB}, we define $G$-type spaces of ultradifferentiable functions $G^{\alpha}_{\alpha}(\RR^d_+)$ and $g^{\alpha}_{\alpha}(\RR^d_+)$ as follows:

Let $A>0$. By $g^{\alpha,A}_{\alpha,A}(\RR^d_+)$, we denote the
space of all $f\in\SSS(\RR^d_+)$ such that
\begin{eqnarray*}
\sup_{p,k\in\NN^d_0}\frac{\|x^{(p+k)/2} \partial^pf(x)\|_{L^2(\RR^d_+)}}
{A^{|p+k|}k^{(\alpha/2)k}p^{(\alpha/2)p}}<\infty.
\end{eqnarray*}
With the following seminorms
\begin{eqnarray*}
\sigma_{A,j}^{\alpha, \alpha} (f)=\sup_{p,k\in\NN^d_0}\frac{\|x^{(p+k)/2} \partial^pf(x)\|_{L^2(\RR^d_+)}}
{A^{|p+k|}k^{(\alpha/2)k}p^{(\beta/2)p}}+\sup_{\substack{|p|\leq
j\\ |k|\leq j}} \sup_{x\in\RR^d_+}|x^k D^pf(x)|,\,\, j\in\NN_0,
\end{eqnarray*}
it becomes an $(F)$-space.

The spaces $G^{\alpha}_{\alpha}(\RR^d_+)$ and
$g^{\alpha}_{\alpha}(\RR^d_+)$ are then defined as union, respectively intersection of
$g^{\alpha,A}_{\alpha,A}(\RR^d_+)$ with respect to $A$:
\begin{equation}
G^{\alpha}_{\alpha}(\RR^d_+)=
\bigcup_{A>0}
g^{\alpha,A}_{\alpha,A}(\RR^d_+), \qquad
g^{\alpha}_{\alpha}(\RR^d_+)=
\bigcap_{A>0}
g^{\alpha,A}_{\alpha,A}(\RR^d_+).
\end{equation}
Then $G^{\alpha}_{\alpha}(\RR^d_+)$ and
$g^{\alpha}_{\alpha}(\RR^d_+)$ are endowed with inductive, respectively projective limit topology.

We note that  $G^{\alpha}_{\alpha}(\RR^d_+)$ is non-trivial when $ \alpha \geq 1$, and
$g^{\alpha}_{\alpha}(\RR^d_+)$ is non-trivial when $ \alpha > 1$.

\par

Next, we rephrase and summarize results from \cite[Theorem 3.6]{D1} and \cite[Theorem 5.7]{SSB} as follows:

\begin{thm} \label{D Th3.6}
Let  $f\in\SSS(\RR^d_+)$.
Then the following is true.
\begin{enumerate}
\item[$(i)$] Let $ \alpha \geq 1$. If $f\in G_\alpha^\alpha(\mathbb{R}^d_+)$, and
$\{a_n\}_{n\in\mathbb{N}_0^d}$ is the sequence of its Laguerre coefficients,
then    $\{a_n\}_{n\in\mathbb{N}_0^d} \in \ell _{\alpha}(\NN_0^d)$.

Conversely if $\{a_n\}_{n\in\mathbb{N}_0^d} \in \ell _{\alpha}(\NN_0^d)$,
then there exists $f\in G_\alpha^\alpha(\mathbb{R}^d_+)$
such that $\{a_n\}_{n\in\mathbb{N}_0^d}$ is the sequence of its Laguerre coefficients.
\item[$(ii)$] Let $ \alpha > 1$. If $f\in g_\alpha^\alpha(\mathbb{R}^d_+)$, and
$\{a_n\}_{n\in\mathbb{N}_0^d}$ is the sequence of its Laguerre coefficients,
then
  $\{a_n\}_{n\in\mathbb{N}_0^d} \in \ell _{0,\alpha}(\NN_0^d)$.
\end{enumerate}

Conversely if $\{a_n\}_{n\in\mathbb{N}_0^d} \in \ell _{0,\alpha}(\NN_0^d)$,
then there exists  $f\in g_\alpha^\alpha(\mathbb{R}^d_+)$
such that $\{a_n\}_{n\in\mathbb{N}_0^d}$ is the sequence of its Laguerre coefficients.
\end{thm}

By Theorems \ref{D Th3.6} and \ref{glavna} we have the following:

\begin{cor} \label{gtype=pil}
 $G$-type spaces coincide with Pilipovi\'c spaces whenever the former are non-trivial. More precisely,
$\mathbf G^\alpha _\alpha(\RR^d_+) = G_\alpha^\alpha(\RR^d_+)$, when $\alpha\geq 1$ and
$\mathbf g ^\alpha _\alpha (\RR^d_+)= g_\alpha^\alpha(\RR^d_+)$, when $\alpha>1$.
\end{cor}

Therefore, $G$-type spaces can be also defined by the iterates of the Laguerre operator, which gives a new insight
into their structure.

It should be noticed that $\mathbf{g}_1^1(\RR^d_+)$ is of a particular interest, since
it is the largest Pilipovi\'c space on orthant which is not a $G$-type space.

\section*{Acknowledgement}
This research was supported by the Science Fund of the Republic of
Serbia, $\#$GRANT No. 2727, {\it Global and local analysis of operators and
	distributions} - GOALS. The first author is supported by the Ministry of Science, Technological Development and Innovation of the Republic of Serbia Grant No. 451-03-65/2024-03/ 200169. The second author is supported by the project F10 financed by the Serbian Academy of Sciences and Arts. The third author is supported by the Ministry of Science, Technological Development and Innovation of the Republic of Serbia (Grants No. 451--03--66/2024--03/200125 $\&$ 451--03--65/2024--03/200125).

\end{document}